\documentclass{article}

\usepackage{amsmath, amsfonts, amsthm} 
\usepackage{algorithm,algorithmic}
\usepackage{graphicx}

\newcommand{\set}[1]{ \left\{ #1 \right\} }
\newcommand{\cH}{\mathcal{H}} % calligraph H for Hilbert spaces
\DeclareMathOperator{\diag}{diag}

\newcommand{\beas}{\begin{eqnarray*}}
\newcommand{\enas}{\end{eqnarray*}}

\newcommand{\R}{\mathbb{R}}
\newcommand{\kernel}{K}
\def\[{\begin{equation}}
\def\]{\end{equation}}

\newcommand{\RR}{\mathbb{R}}
\newcommand{\N}{\mathbb{N}}

\usepackage{graphicx}

\newtheorem{theorem}{Theorem}

\usepackage{authblk}
\title{Greedy Kernel Methods for Center Manifold Approximation}
\author[2]{B. Hamzi \thanks{boumediene.hamzi@gmail.com}}
\affil[2]{Department of Mathematics, AlFaisal University, Riyadh, KSA}

\author[1]{B. Haasdonk \thanks{bernard.haasdonk@mathematik.uni-stuttgart.de}}
\author[1]{G. Santin \thanks{gabriele.santin@mathematik.uni-stuttgart.de, 
orcid.org/0000-0001-6959-1070}}
\author[1]{D. Wittwar \thanks{dominik.wittwar@mathematik.uni-stuttgart.de}}
\affil[1]{Institute for Applied Analysis and Numerical Simulation, University of Stuttgart, Germany}

\begin{document}

\maketitle

\abstract{
For certain dynamical systems it is possible to significantly simplify the study of stability by means of the center manifold theory.
This theory allows to isolate the complicated asymptotic behavior of the system close to a non-hyperbolic equilibrium point, and to obtain 
meaningful predictions of its behavior by analyzing a reduced dimensional problem.
Since the manifold is usually not known, approximation methods are of great interest to obtain qualitative estimates. 
In this work, we use a data-based greedy kernel method to construct a suitable approximation of the manifold close to the equilibrium. The data are collected by repeated 
numerical simulation of the full system by means of a high-accuracy solver, which generates sets of discrete trajectories that are then used to construct a surrogate 
model of the manifold. 
The method is tested on different examples which show promising performance and good accuracy.
}

% BEGIN OF INTRODUCTION

\section{Introduction}
Center manifold theory plays an important role in the study of the
stability of dynamical systems when the equilibrium point is not
hyperbolic.  It  isolates the complicated asymptotic behavior by locating the center manifold which is an
invariant manifold tangent to the subspace spanned by the
eigenspace of eigenvalues on the imaginary axis. Then, the dynamics of the original system
 will be essentially determined by the  restriction of this dynamics on the center manifold 
since the local dynamic behavior ``transverse'' to this invariant
manifold is relatively simple as it corresponds to the flows in the local stable (and unstable) manifolds. In practice, one
does not compute the center manifold and its dynamics exactly
since this requires the resolution of a quasilinear partial
differential  equation which is not easily solvable. In most cases
of interest, an approximation of degree two or three of the
solution is sufficient. Then, the reduced dynamics on
the center manifold can be determined, its stability can be studied
and then conclusions about the stability of the
original system can be obtained \cite{pliss, sositaisvili,
kelley, carr81, henry}.

In this article, we use greedy kernel methods to construct a data-based approximation of the center manifold. The present work is a preliminary study that 
is intended to introduce our concept and algorithm, and to test it on some examples.

% END OF INTRODUCTION

% BEGIN OF BACKGROUND

\section{Background}

Consider a dynamical system
\begin{equation}\dot{x}=f(x)=Fx+\bar{f}(x) \label{eqn:fors}
\end{equation}
of large dimension $n$, where $f : D \rightarrow \RR^n$ is a continuously differentiable function over the domain $D \subset \R^n$ such that $0 \in D$. Let 
$F=\frac{\partial\,f}{\partial\,x}(x)|_{x=0}$. Suppose $x=0$ is an equilibrium, i.e. $f(0)=0$, and denote as $\sigma_{\R}(F)$ the sequence of real parts of the 
eigenvalues of $F$. 
A classical stability result states that if $F$ has all its eigenvalues with negative real parts, i.e., $\sigma_{\R}(F)\subset \R_{<0}$, then the origin is 
asymptotically stable; and if $F$ has some eigenvalues with positive real parts, then the origin is unstable. If  $\sigma_{\R}(F)\subset \R_{\leq0}$, the linearization 
fails to determine the stability properties of the origin.

  After a linear change of 
coordinates, we have
\begin{eqnarray} \label{full_order_system}
\dot{x}&=& F_1 x + \bar{f}_1(x,y) \label{fors1} \\
\dot{y}&=& F_2 y + \bar{f}_2(x,y) \label{fors2}
\end{eqnarray}
where $F_1 \in \RR^{d \times d}$ is such that $\sigma_{\R}(F_1)=\{0\}$ and $F_2 \in \RR^{m \times m}$ with $m := n-d$ is such that 
$\sigma_{\R}(F_2)\subset\R_{<0}$.   The functions $\bar{f}_1 : \R^d \times \R^m \rightarrow \R^d$ and $\bar{f}_2 : \R^d \times \R^m \rightarrow \R^m$ are continuously 
differentiable. 

Intuitively we expect the 
stability of the
equilibrium to only depend on the nonlinear terms $\bar{f}_1(x,y) $.  The center manifold theorem correctly formalizes
this intuition.

A center manifold is an invariant manifold\footnote{A differentiable manifold $\cal{M}$ is said to be invariant under the flow of a vector field $X$ if for $x \in \cal{M}$, $F_t(x) \in \cal{M}$ for small $t > 0$, where $F_t(x)$ is the flow of
$X$.}, $ y=h(x)$,  for (\ref{fors1})--(\ref{fors2}), such that $h$ is smooth and
\[\label{eq: center manifold properties}h(0)=0; \quad Dh(0)=0. \]
%\gs{@Boumediene: can you please provide a reference for this theorem?}
\begin{theorem} \cite{carr81}
 If $f_1$ and $f_2$ are twice continuously differentiable and are such that
 \[\bar{f}_i(0,0)=0; \frac{\partial f_i}{\partial x}(0,0)=0; \frac{\partial f_i}{\partial y}(0,0)=0,  \]
 for $i=1,2$, and the eigenvalues of $F_1$ have zero real parts, and all the eigenvalues of $F_2$ have negative real parts, then there exists a neighbourhood $\Omega 
\subset \R^d$ of
the origin $0 \in \R^d$, such that $ x_2=h(x_1)$ is a center manifold for  (\ref{fors1})--(\ref{fors2}).
\end{theorem}

 Since  
\begin{eqnarray*}
 \dot{x}&= F_1 x + \bar{f}_1(x,y), \\
\dot{y} &= F_2 y + \bar{f}_2(x,y)
\end{eqnarray*}
 and $x_2=h(x_1)$, we deduce that $h $ satisfies the PDE

\begin{equation}\label{cm_pde} F_2 h(x) + \bar{f}_2(x,h(x))=D 
h(x)
\left(F_1 x +\bar{f}_1(x,h(x))\right).\end{equation}

 The center manifold theorem ensures that there are smooth solutions 
to this PDE. It also allows to deduce the stability of the origin of the full order system (\ref{fors1})--(\ref{fors2}) from the stability of the origin of a reduced 
order system called the \emph{center dynamics}.
 
%\gs{@Boumediene: could you please add a reference for this Theorem?} 

\begin{theorem} \cite{carr81} \label{cm_theorem} (Center Manifold Theorem) The equilibria $x=0,y=0$ of the original dynamics is locally 
asymptotically stable (resp. unstable)
iff the equilibria $x=0$ of the center dynamics (dynamics on the center manifold)
\begin{equation}\label{reduced_system} \dot{x}= F_1 x +\bar{f}_1(x,h(x)),\end{equation}
is locally  asymptotically stable (resp. unstable).
\end{theorem}

  After solving
the PDE (\ref{cm_pde}),  the problem of analyzing the stability properties of the system (\ref{full_order_system})  reduces to analyzing the nonlinear 
stability of a lower dimensional system (\ref{reduced_system}).  

But the PDE (\ref{cm_pde}) need not be solved exactly, since frequently it suffices to compute the low degree 
terms of the Taylor series expansion of $h$ around
$x=0$, i.e., 
 
$$h(x)=h^{[1]} x+h^{[2]}(x)+h^{[3]}(x)+\ldots$$
where $(\cdot)^{[k]}$ is the degree $k$ part of the Taylor series. The PDE can then be rewritten as a set of algebraic equations
 \begin{itemize}
\item  $F_2 h^{[1]}=h^{[1]} F_1$
\item  $ \scriptstyle F_2 h^{[2]}(x) + \bar{f}^{[2]}_2(x,h^{[1]}(x))=\scriptstyle \frac{\partial 
h^{[2]}}{\partial x}(x)\left(F_1 x_1 +\bar{f}^{[2]}_1(x,h^{[1]}(x))\right)$
\item 
  $ \scriptstyle F_2 h^{[3]}(x) + 
\left(\bar{f}_2(x,h^{[2]}(x))\right)^{[2]}=\frac{\partial h^{[2]}}{\partial
x}(x)
\left(F_1 x +\left(\bar{f}_1(x,h^{[2]}(x))\right)^{[2]}\right)$
\item
etc.
\end{itemize}

The Taylor expansion of $h$ to degree $d-1$ determines the 
center manifold dynamics
to degree $d$,
\begin{eqnarray*} \dot{x}&=&F_1 x
\\& +&\left(\bar{f}_1(x,h^{[1]} x)\right)^{[2]}
\\&+&\left(\bar{f}_1(x,h^{[1]} x+h^{[2]}(x))\right)^{[3]}
\\&+& \ldots
\\&+&\left(\bar{f}_1(x,h^{[1]} x+h^{[2]}(x)+\ldots+h^{[d-1]}(x))\right)^{[d]},
\end{eqnarray*}
and this may be enough to determine its local asymptotic stability. This methodology is valid for parameterized dynamical systems and is used to study the stability of 
dynamical systems with bifurcations. 

Our goal in this paper is to find a data-based approximation of the center manifold in view of a data-based version of the center manifold theorem.

%\section{Setting}
%
%We consider autonomous systems of the form
%\begin{equation}
% \dot{u} = f(u) \label{eq: autonomous system}
%\end{equation}
%where $f : D \rightarrow \RR^n$ is a continuously differentiable function over the domain $D \subset \R^n$ such that $0 \in D$.
%We further assume that $u = 0$ is an equilibrium point of \eqref{eq: autonomous system}, and that the derivative 
%$Df : D \rightarrow \R^{n \times n}$ of $f$ evaluated at the origin has $d$ eigenvalues with zero real part and $m := n - d$ 
%eigenvalues with negative real part. Due to this we may assume that after an appropriate coordinate transformation the system
%\eqref{eq: autonomous system} takes the form
%\begin{equation} \label{eq: transformed system}
% \begin{split}
%  \dot{x} & = f_1(x,y) \\
%  \dot{y} & = f_2(x,y) \\
% \end{split}
%\end{equation}
%with $f_1 : \R^d \times \R^m \rightarrow \R^d$ and $f_2 : \R^d \times \R^m \rightarrow \R^m$ continuously differentiable. In
%\eqref{eq: transformed system} the $d$ components in $x$ correspond to the eigenvalues with zero real part and the $m$ components
%in $y$ correspond to the eigenvalues with negative real part. By this we have for a suitable neighbourhood $\Omega \subset \R^d$ of
%the origin $0 \in \R^d$ that there exists a center manifold $h : \Omega \rightarrow \R^m$, i.e. $h$ is continuously differentiable
%and it holds
%\begin{equation}
% h(0) = 0,  \quad Dh(0) = 0 \quad \text{and} \quad y = h(x). \label{eq: center manifold properties}
%\end{equation}

% END OF BACKGROUND

% BEGIN OF KERNEL APPROXIMATION

\section{Kernel Approximation}
We want to build a surrogate model $s_h :\Omega\rightarrow \R^m$ which approximates the center manifold $h$ on a suitable set $\Omega\subset\R^d$, in the 
sense that $s_h(x) \approx h(x)$ for all $x\in\Omega$. This model is constructed in a data-based way, i.e., we assume the knowledge of the map $h$ on a finite 
set of input parameters, or training data. In practice, such values are computed from high-fidelity numerical approximations, which will be discussed in details in the 
following. 

The surrogate is based on kernel approximation, which allows the use of scattered data, i.e., we do not require any grid structure on the set of training data. Moreover, 
since the unknown function $h$ is vector-valued, we employ here matrix-valued kernels. Details on kernel-based approximation can be found e.g. in \cite{Wendland2005}, 
and the extension to the vectorial case is detailed e.g. in \cite{Micchelli2005, 2018arXiv180709111W}. We recall here only that a positive definite matrix-valued kernel 
on $\Omega$ is a function $\kernel : \Omega \times \Omega \rightarrow \R^{m \times m}$ such that $\kernel(x, y) = \kernel(y,x)^T$ for all $x,y \in \Omega$ and 
$[\kernel(x_i, x_j)]_{i,j=1}^N\in\R^{mN\times 
mN}$ is positive semidefinite for any set $\{x_1, \dots, x_N\}\subset\Omega$ of pairwise distinct points, for all $ N \in\N$. Associated to a positive definite kernel 
there is a unique Hilbert 
 space $\mathcal H$ of functions $\Omega\to\R^m$, named native space, where the kernel is reproducing, meaning that $\kernel(\cdot, x) \alpha$ is the Riesz representer 
of the 
directional point evaluation $\delta_x^{\alpha}(f):= \alpha^T f(x)$, for all $\alpha\in\R^m$, $x\in\Omega$. 

We consider here a twice continuously differentiable matrix-valued kernel $k$ on $\Omega$, and we use a specific functional 
formulation for our approximation and a specific cost function, in order to construct a surrogate that is well suited for the particular approximation task. 

In details, the approximant takes the form
\begin{equation*}
 s_h(x) = \sum\limits_{i=1}^{n_1} \kernel(x,x_i^{(1)})\alpha_i + \sum\limits_{j=1}^{n_2} \sum\limits_{i=1}^m \partial_i^{(2)} \kernel(x,x_j^{(2)})\beta_{i,j},
\end{equation*}
with centers $x_i^{(1)} \in X^{(1)} = \set{ x_1^{(1)}, \dots, x_{n_1}^{(1)} }$, $x_j^{(2)} \in X^{(2)} = \set{ x_1^{(2)}, \dots, x_{n_2}^{(2)} }$ and coefficient vectors
$\alpha_i, \beta_{i,j} \in \R^m$. Here the superscript $\partial^{(2)}$ denotes that the derivative with regards to the second kernel
component is taken.

We assume to have sufficiently many data $X_{N^{\ast}} = \{x_1, \dots,x_{N^{\ast}} \}$ and $Y_{N^{\ast}} = \{ y_1, \dots, y_{N^{\ast}} \}$ which, for example, are 
generated by running
a numerical scheme to compute discrete trajectories for different initial values $(x_0,y_0)$. For this step, we need to assume that the variable splitting 
\eqref{fors1}--\eqref{fors2} is known in advance.
Note that this is not a severe restriction, as for a general ODE
\eqref{eqn:fors} the required state transformation can be determined by
eigenvalue decomposition of $F$.

Observe that we do not know if a data pair $(x_i,y_i)$ lies on the center manifold, i.e. if $y_i = h(x_i)$ holds. We only know that the data converges 
asymptotically to the center manifold as $x_i \rightarrow 0$. Thus, an interpolation-based surrogate which interpolates the data
on a given subset $X \subset X_{N^{\ast}}$ seems ill-suited for our purposes. We consider instead another set of conditions to define the approximant. First, we still 
require the 
conditions in \eqref{eq: center manifold properties} to be satisfied by our approximation. Moreover, for the given subsets $X = \{ x_1 , \dots,x_N \}$ and $Y = 
\{ y_1, \dots, y_N \}$, we compute our approximant by minimizing the following functional $J : \cH \rightarrow \R$ under the constraint $s(0) = 0, Ds(0) = 0$:

\begin{equation}
 J(s) := \| s \|_{\cH}^2 + \sum\limits_{i=1}^N (s(x_i) - y_i)^T \omega_i (s(x_i) - y_i). \label{eq: Functional}
\end{equation}
Here $\omega_i \in \R^{m \times m}$ is a positive definite weight matrix. It can be shown that \eqref{eq: Functional} has a unique minimizer $s_h$ (see \cite{WSH18}).
In particular $s_h$ and its derivative $Ds_h$ have the form
\begin{align}
\label{eq:surrogate} s_h(x) & = \sum\limits_{i=1}^{N+1} \kernel(x,x_i)\alpha_i + \sum\limits_{i=1}^m \partial_i^{(2)}\kernel(x,0)\beta_{i}, \\
\nonumber Ds_h(x) & =  \sum\limits_{i=1}^{N+1} D^{(1)}\kernel(x,x_i)\alpha_i + \sum\limits_{i=1}^m D^{(1)}\partial_i^{(2)}\kernel(x,0)\beta_{i},
\end{align}
where we set $x_{N+1} := 0$. The coefficient vectors $\alpha_i, \beta_i$ can be computed by solving the system
\begin{equation}\label{eq:thelinearsystem}
 \begin{pmatrix}
  A + W & B \\
  B^T & C 
 \end{pmatrix}
 \begin{pmatrix}
  \boldsymbol{\alpha} \\
  \boldsymbol{\beta}
 \end{pmatrix}
  =
 \begin{pmatrix}
  \boldsymbol{f}_1 \\
  \boldsymbol{f}_2
 \end{pmatrix},
\end{equation}
with
\begin{align*}
 A & := \left( \kernel(x_i,x_j) \right)_{i,j} \in \R^{m(N+1) \times m(N+1)}, \\
 W & := \diag\left( \omega_1^{-1},\dots, \omega_N^{-1},0 \right) \in \R^{m(N+1) \times m(N+1)}, \\
 B & := \left( \partial_{j}^{(2)}\kernel(x_i,0) \right)_{i,j} \in \R^{m(N+1) \times m^2}, \\
 C & := \left( \partial_{i}^{(1)}\partial_{j}^{(2)}k(0,0) \right)_{i,j} \in \R^{m^2 \times m^2}, \\
 \boldsymbol{f}_1 & := (y_1^T, \dots,y_n^T,0)^T \in \R^{m(N+1)}, \\
 \boldsymbol{f}_2 & := 0 \in \R^{m^2 \times m}.
 \end{align*}
The weight matrices $\omega_i$ can either be chosen manually, or a regularizing function $r : \Omega \rightarrow \R^{m \times m}$ can be prescribed
such that $\omega_i = r(x_i)$ is symmetric and positive definite. In our numerical examples in section \ref{sec: Numerical examples} we chose a constant regularization 
function, i.e.
\begin{equation*}
 \omega_i = r(x_i) = \lambda I_m
\end{equation*}
for some $\lambda > 0$. However, one might consider a more general approach, where the weight increases as the data tends to the origin,
i.e. $\omega_i \succeq \omega_j$ if $ \|x_i\| \leq \|x_j\|$.

\subsection{Greedy approximation}
If the technique of the previous section is used as it is, the surrogate \eqref{eq:surrogate} is given by an expansion with ${N^{\ast}}$ terms, where ${N^{\ast}}$ is the 
number of 
points in the training set. Therefore, the model evaluation can be not efficient enough if the model is built using a too large dataset. Furthermore, the 
computation of the coefficients in \eqref{eq:surrogate} requires the solution of the linear system \eqref{eq:thelinearsystem}, whose size again scales with the size of 
the training set, and which can be severely ill-conditioned for not well-placed points.

To mitigate both problems, we employ an algorithm that aims at selecting small subsets $X_N$, $Y_N$ of points such that the surrogate computed with 
these sets is a sufficiently good approximation of the one which uses the full sets. The algorithm selects the points in a greedy way, i.e., one point at a time is 
selected and added to the current training set. In this way, it is possible to identify a good set without the need to solve a nearly infeasible combinatorial problem. 

The selection is performed using the $P$-greedy method of \cite{DeMarchi2005} applied to the kernel $K$, such that the set of points is selected before the computation 
of the surrogate.
The number of points, and therefore the expansion size and evaluation speed,
is depending on the prescribed target accuracy $\varepsilon_{tol}>0$. 
For details on the method implementation and its convergence
properties we refer to \cite{SH16b}.

\section{Numerical Examples} \label{sec: Numerical examples}
We test now our method on three different examples. In each of them, we
specify the setting and the parameters used to build the surrogate. We
compare the surrogates with the the true manifold and we compute the
pointwise residual
\begin{equation*}
 r_i(x) = Ds_i(x)f_1(x,s_i(x)) - f_2(x,s_i(x)), \quad i = 1,2,
\end{equation*}
which measures the error when the surrogate is used as a replacement
of the center manifold.

In all the three examples, the greedy algorithm is used to select a suitable subset of the points, and in all cases the procedure is stopped with a prescribed 
$\varepsilon_{tol}$ . In the first two examples we set $\varepsilon_{tol}:=10^{-15}$, while $\varepsilon_{tol}:=10^{-10}$ in the last one.

\subsection{Example 1}

We consider the $2$-dimensional system 
\begin{equation} \label{eq: Example 1 system}
 \begin{split}
  \dot{x} & = f_1(x,y) = xy \\
  \dot{y} & = f_2(x,y) = -y + x^2.  
 \end{split}
\end{equation}
We generate the training data by solving \eqref{eq: Example 1 system} with an implicit Euler scheme for initial time $t_0 = 0$, final time
$T = 1000$ and with the time step $\Delta t = 0.1$. We initiate the numerical procedure with initial values $(x_0,y_0) \in \set{ \pm 0.8 } \times \set{ \pm 0.8}$
and store the resulting data pairs in $X$ and $Y$ after discarding all data whose $x$-values are not contained in the neighborhood $[-0.1,0.1]$ 
which results in $N^\ast = 38248$ data pairs. 

\begin{figure}[ht]
\begin{center}
\begin{tabular}{c}
\includegraphics[width=0.8\textwidth]{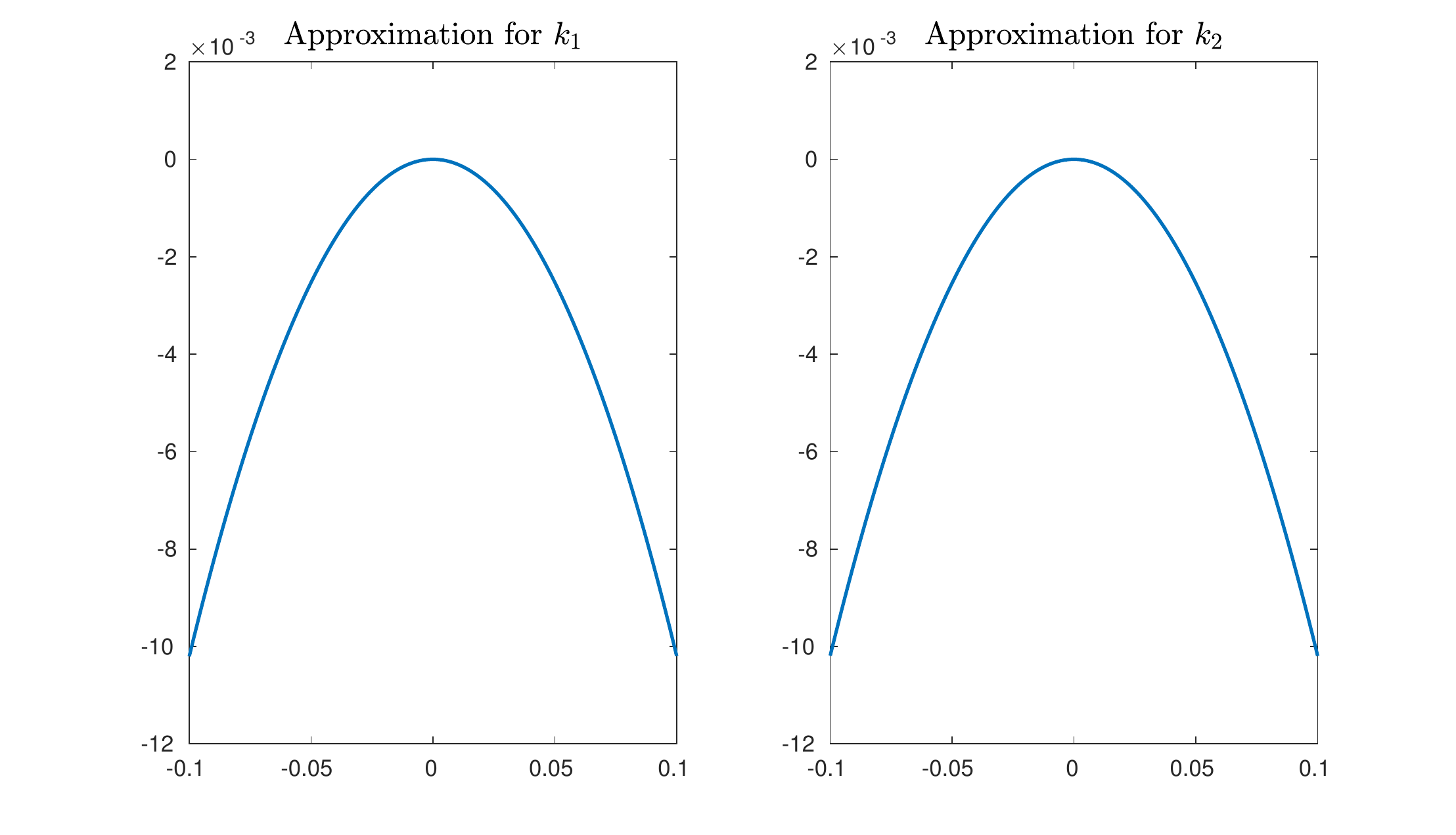}\\
\includegraphics[width=0.8\textwidth]{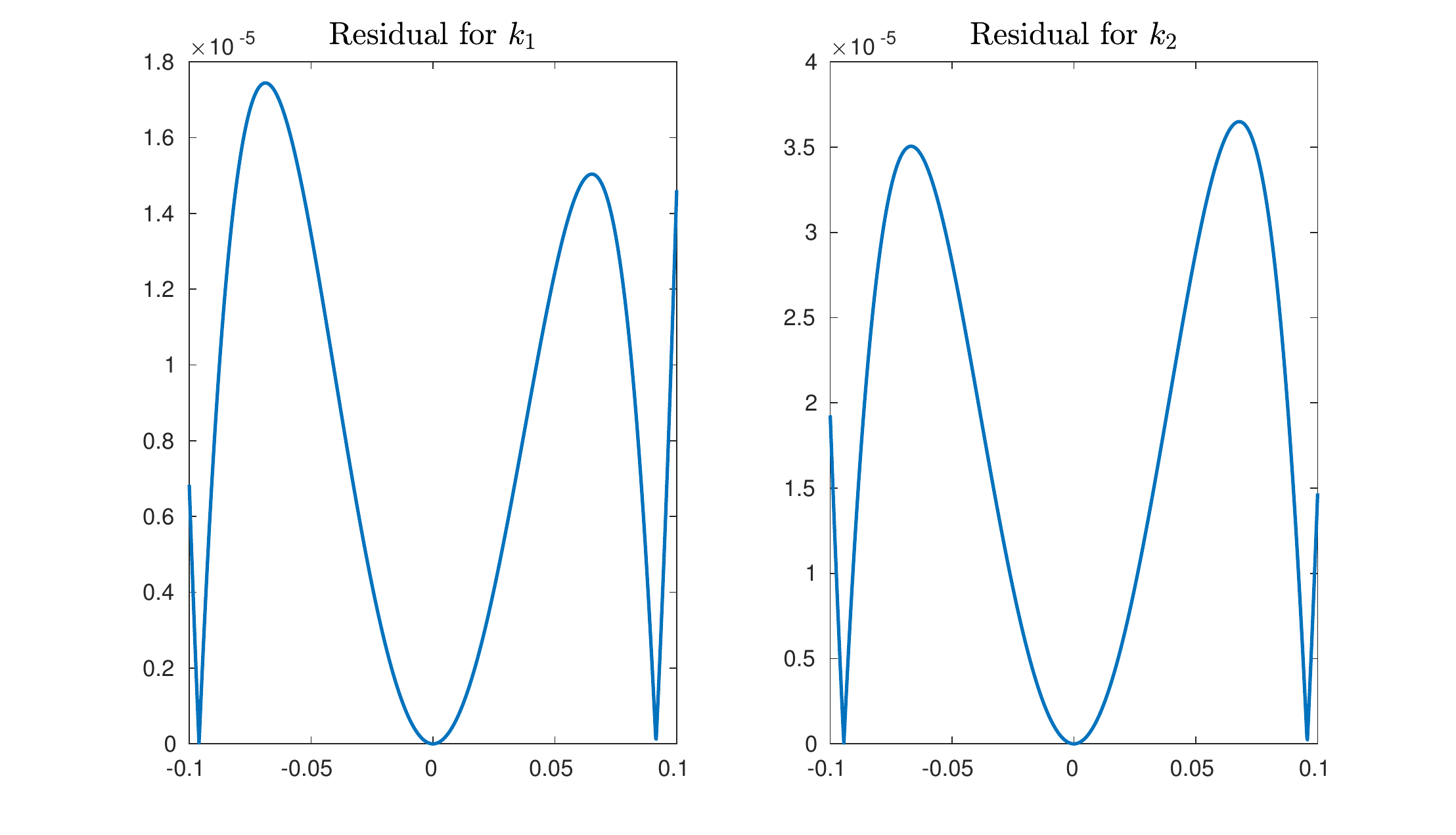}
\end{tabular}
\end{center}
\caption{Approximation of the center manifold and residuals for Example 1 using kernels $k_1$ and $k_2$.}
\label{figure:1}
\end{figure}

We run the greedy algorithm for the kernels $k_1(x,y) := \left(1 + xy/2\right)^4$ and $ k_2(x,y) = e^{-(x-y)^2/2}$. This results in the sets $X_1$ and $X_2$ which 
contain 
$14$ and $6$ points, respectively. The corresponding approximations $s_1$ and $s_2$ for the constant regularization function $r \equiv 10^{-10}$ are plotted in Figure 
\ref{figure:1}, left. The center manifold approximations are plotted over the domain $[-0.1,0.1]$. The pointwise residual is depicted in Fig.\ \ref{figure:1}, right.

% Figures Example 1

\subsection{Example 2}

We consider the $2$-dimensional system
\begin{equation} \label{eq: Example 2 system}
 \begin{split}
  \dot{x} & = f_1(x,y) = -xy \\
  \dot{y} & = f_2(x,y) = x^2 -y - 2y^2.  
 \end{split}
\end{equation}

\begin{figure}[ht]
\begin{center}
\begin{tabular}{c}
\includegraphics[width=0.8\textwidth]{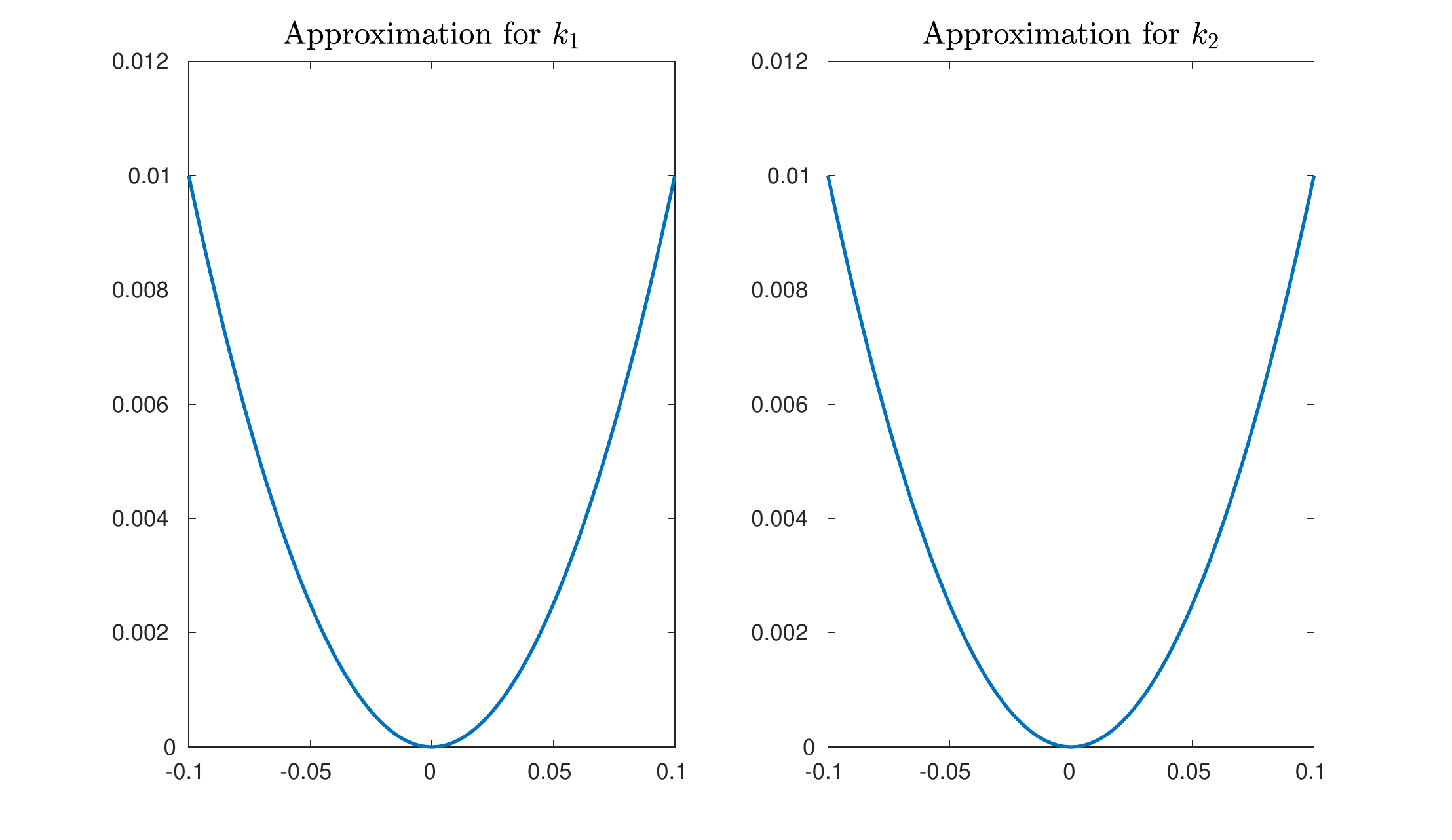}\\
\includegraphics[width=0.8\textwidth]{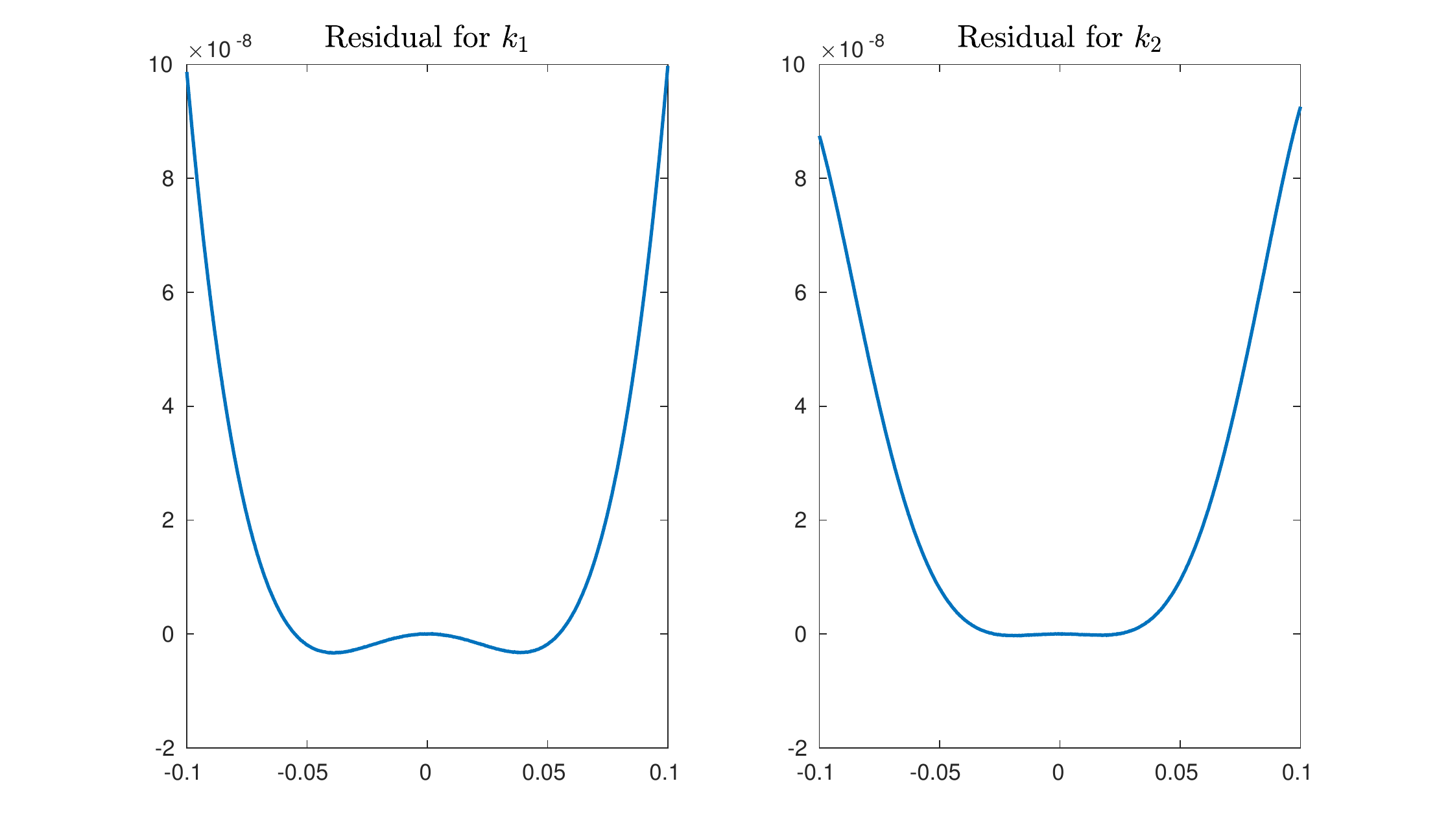}
\end{tabular}
\end{center}
\caption{Approximation of the center manifold and residuals for Example 2 using kernels $k_1$ and $k_2$.}
\label{figure:2}
\end{figure}

The training data is generated the same way as in Example 1. We again use the kernels $k_1$ and $k_2$. The greedy algorithm gives sets $X_1$ and $X_2$ of size $12$ 
and $6$, respectively. The evaluation of the approximations $s_1$ and $s_2$ over the neighborhood $[-0.1,0.1]$ can be seen on the left in Figure \ref{figure:2}, while 
the 
respective pointwise residuals are plotted in the right part of Figure \ref{figure:2}.

% Figures Example 2

\subsection{Example 3}
We consider the $(2+1)$-dimensional system
\begin{equation} \label{eq: Example 3 system}
 \begin{split}
  \dot{x} & = f_1(x,y) = \begin{pmatrix}
               -x_2 + x_1y \\ x_1 + x_2y 
              \end{pmatrix}\\
    \dot{y} & = f_2(x,y) = -y -x_1^2 - x_2^2 + y^2.
 \end{split}
\end{equation}

\begin{figure}[ht]
\begin{center}
\begin{tabular}{c}
\includegraphics[width=0.8\textwidth]{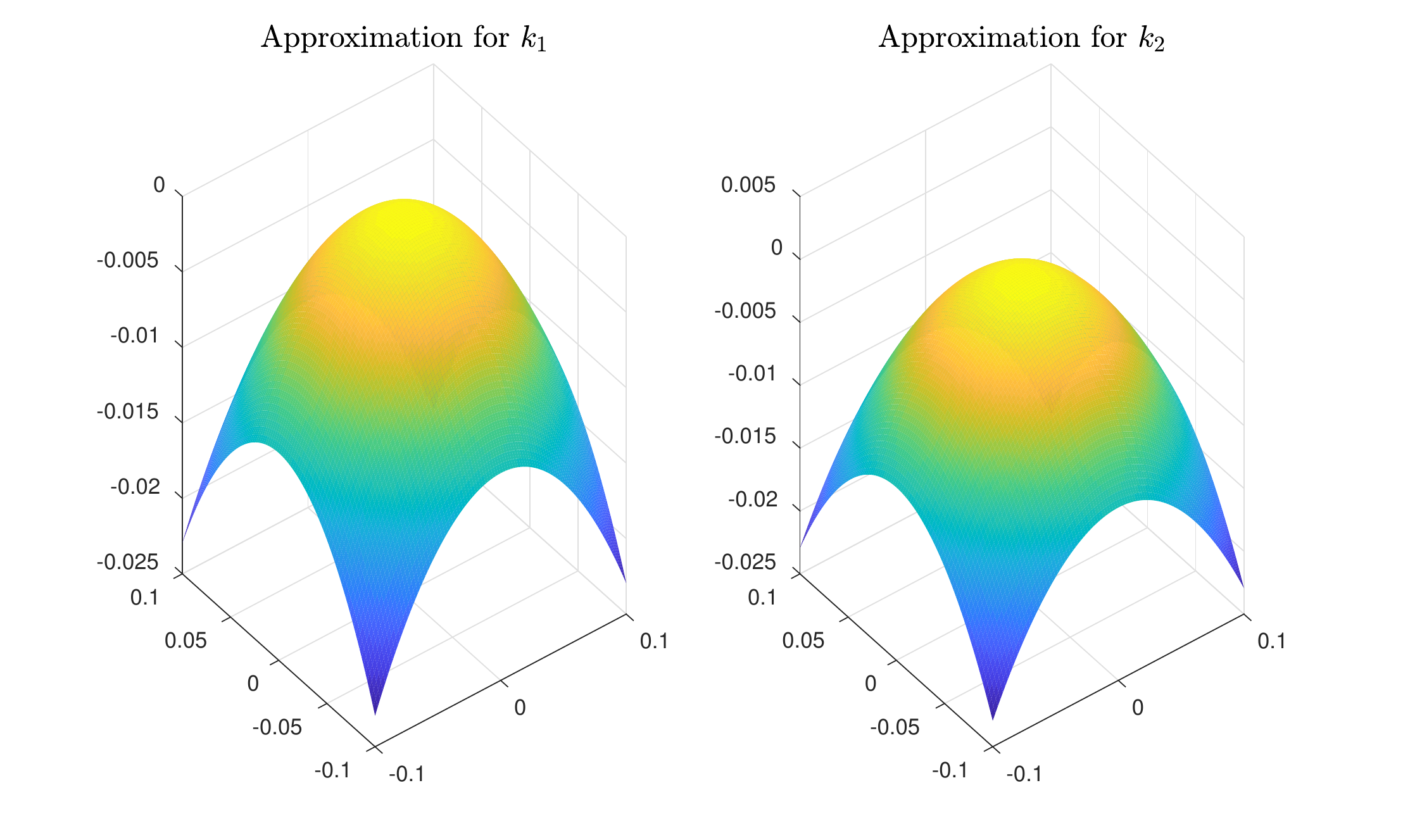}\\
\includegraphics[width=0.8\textwidth]{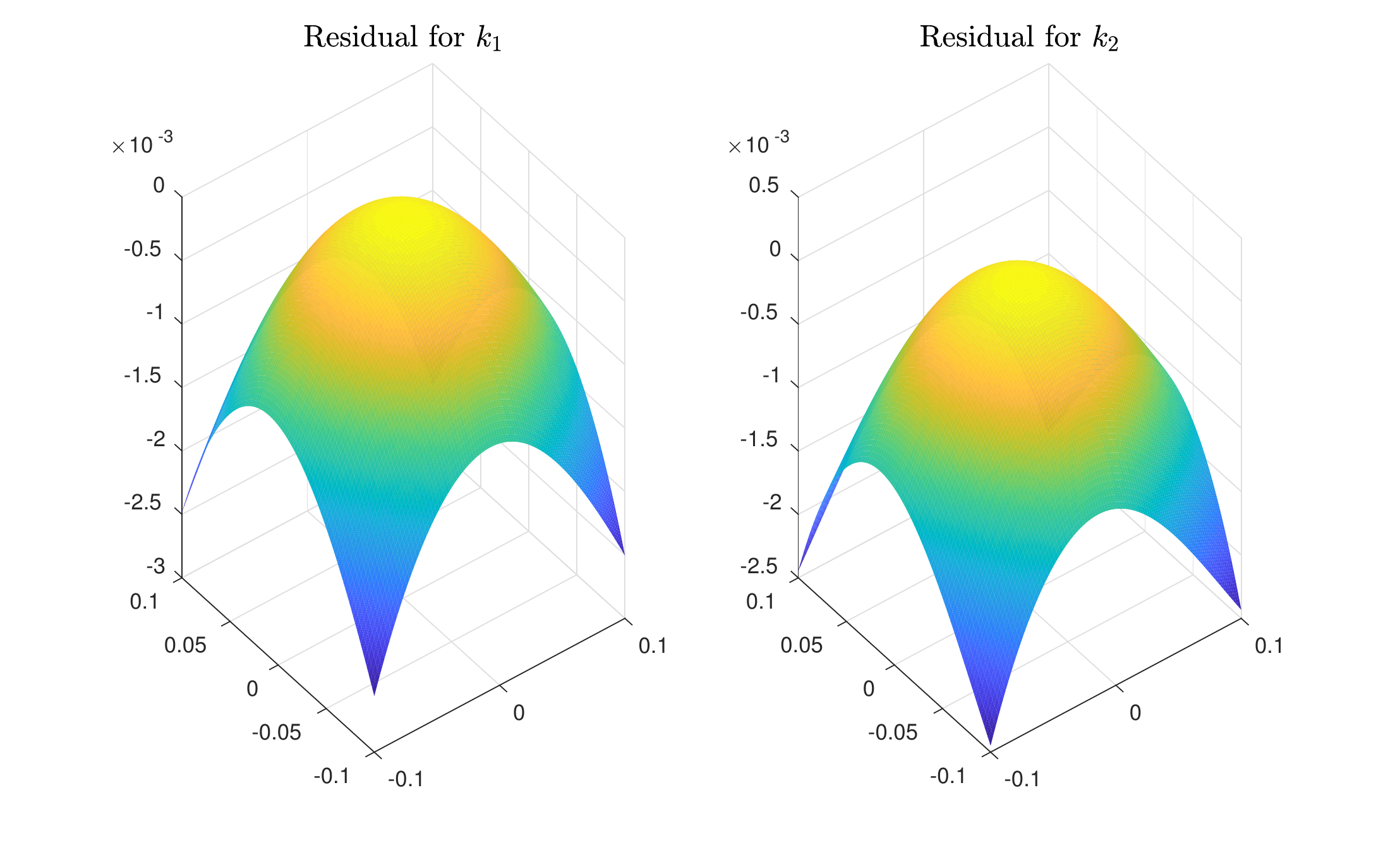}
\end{tabular}
\end{center}
\caption{Approximation of the center manifold and residuals for Example 3 using kernels $k_1$ and $k_2$.}
\label{figure:3}
\end{figure}

We generate the training data in a similar fashion as before. We again use the implicit Euler scheme with start time $t = 0$, final time $T = 1000$ and with
time step $\Delta t = 0.1$. The Euler method is performed for initial data $(x_0,y_0) \in \set{ \pm 0.8 }^3$ and the resulting trajectories are stored in $X$ and $Y$,
where only data with $x \in [-0.1,0.1]^2$ was considered; this leads to $N^\ast = 78796$ data pairs. We use the kernels $k_1(x,y) = (1 + x^Ty/2)^4$ and $k_2(x,y) = e^{- 
\|x - 
y\|_2^2/2}$, and the greedy-selected sets have the size $21$ (for $k_1$) and $25$ (for $k_2$), respectively. The approximations $s_1$,$s_2$ and
their corresponding residuals $r_1$ and $r_2$ computed over the domain $[-0.1,0.1]^2$. The results can be seen in Figure \ref{figure:3}.

% Figures Example 3

We remark that in all the three experiments both kernels give comparable results in terms of error magnitude, and they both provide a good approximation of the manifold.

\section{Conclusions}
In this paper we introduced a novel algorithm to approximate the
center manifold of a given ODE using a data-based surrogate. 

This algorithm computes an approximation of the manifold from a set of numerical trajectories with different initial data. It is based on kernel 
methods, which allow the use of the scattered data generated by these simulations as training points. Moreover, an application-specific ansatz and cost function have 
been employed in order to enforce suitable properties on the surrogate.

Several numerical experiments suggested that the present method can reach a significant accuracy, and that it has the potential to be used as an effective model 
reduction technique.
It seems promising to apply this approach to high dimensional systems as the
approximation technique straightforwardly can be extended and is less prone to the curse of dimensionality than grid-based approximation techniques.
An interesting extension would consist of 
determining the
decomposition \eqref{fors1}--\eqref{fors2} in a data-based fashion
by suitable processing of the trajectory data.

\section{Acknowledgements} 
The first, third, and fourth authors would like to thank the German Research Foundation (DFG) for support within the Cluster of Excellence in
Simulation Technology (EXC 310/2) at the University of Stuttgart.

\bibliographystyle{abbrv}
\bibliography{biblio}

\end{document}